\documentclass[a4paper,12pt]{amsart}
\usepackage{amssymb}
\usepackage{amsmath}
\usepackage{amstext}
\usepackage{amsgen}
\usepackage{amsbsy}
\usepackage{amsopn}

\newtheorem{e-proposition}[theorem]{Proposition}

\newtheorem{e-definition}[theorem]{Definition\rm}

\def\FM{{\mathbb{F}}}

\def\QM{{\mathbb{Q}}}

\def\ZM{{\mathbb{Z}}}

\def\Bb{{\mathbf B}}
\def\Gb{{\mathbf G}}
\def\Lb{{\mathbf L}}
\def\Pb{{\mathbf P}}
\def\Tb{{\mathbf T}}
\def\Ub{{\mathbf U}}
\def\Vb{{\mathbf V}}
\def\Xb{{\mathbf X}}

\def\BCB{{\boldsymbol{\mathcal{B}}}}
\def\OCB{{\boldsymbol{\mathcal{O}}}}
\def\vdo{{\dot{v}}}
\def\wdo{{\dot{w}}}
\def\xdo{{\dot{x}}}

\def\longto{\longrightarrow}

\def\fonctio#1#2#3#4{\begin{array}{ccc}
{#1} & \longto & {#2} \\
{#3} & \longmapsto & {#4} 
\end{array}}

\def\fin{~$\scriptstyle \blacksquare$}

\def\lexp#1#2{\kern\scriptspace\vphantom{#2}^{#1}\kern-\scriptspace#2}

\def\gfq{\FM_{\! q}}
\def\cyclic{\hspace{0.2em}\smash{
     {\mathop{\longleftrightarrow}\limits^{F}}\hspace{0.2em}}}
\def\qlb{\overline{\QM}_{\! \ell}}

\begin{document}

\baselineskip=16pt
\title{Affineness of Deligne-Lusztig varieties \\ for minimal length elements}

\author{C\'edric Bonnaf\'e},
\address{Laboratoire de Math\'ematiques de Besan\c{c}on
(CNRS - UMR 6623) \\
Universit\'e de Franche-Comt\'e, 16 Route de Gray,
25030 Besan\c{c}on Cedex, France}
\makeatletter
\email{cedric.bonnafe@math.univ-fcomte.fr}

\author{Rapha\"el Rouquier}
\address{Mathematical Institute, University of Oxford\\
24-29 St Giles', Oxford, OX1 3LB, UK} 
\makeatletter
\email{rouquier@maths.ox.ac.uk}

\maketitle

\pagestyle{myheadings}

\markboth{\sc C.Bonnaf\'e \& R. Rouquier}{\sc 
Affineness of Deligne-Lusztig varieties}

\begin{abstract}
We prove that the Deligne-Lusztig varieties associated to elements 
of the Weyl group which are of minimal length in their twisted class 
are affine. Our proof differs from the proof of He
and Orlik-Rapoport and it is inspired by the case of regular elements, which
correspond to the varieties involved in Brou\'e's conjectures.

\end{abstract}

\section{Introduction}
\label{}

Let $p$ be a prime number, let $\FM$ denote an algebraic closure of the finite 
field with $p$ elements and let $\Gb$ be a connected reductive algebraic 
group over $\FM$. We assume that $\Gb$ is endowed with an isogeny
$F : \Gb \to \Gb$ such that $F^\delta$ is a Frobenius endomorphism with
respect to some $\gfq$-structure 
on $\Gb$ (here, $\delta$ is a non-zero natural number, $q$ is a power of $p$
and $\gfq$ denotes the subfield of $\FM$ with $q$ elements). 

We denote by $\BCB$ the variety of Borel subgroups of $\Gb$ and by
$\BCB\times\BCB=\coprod_{w\in W}\OCB(w)$ the decomposition into orbits for the
diagonal action of $\Gb$. Here, $W$ is the Weyl group of $\Gb$, with set
of simple reflections $S$ corresponding to the orbits of dimension
$1+\dim\BCB$, and the first and last projections define an isomorphism
$\OCB(w)\times_\BCB\OCB(w')\xrightarrow{\sim}\OCB(ww')$ when
$\ell(ww')=\ell(w)+\ell(w')$, where
$\ell : W \to \ZM_{\ge 0}$ is the length function on $W$ associated to $S$.

Given $w \in W$, we define the {\it Deligne-Lusztig variety} 
\cite[Definition 1.4]{delu} associated to $w$ by 
$$\Xb(w)=\Xb_\Gb(w)=\{\Bb \in \BCB~|~(\Bb,F(\Bb))\in \OCB(w)\}.$$

\medskip
By studying a class of ample sheaves on $\Xb(w)$, 
Deligne and Lusztig proved that these varieties are affine 
when 
$q^{1/\delta}$ is larger than the Coxeter number 
of $\Gb$ \cite[Theorem 9.7]{delu}.

\smallskip

They proved more generally that the existence of coweights satisfying certain
inequalities ensures that $\Xb(w)$ is affine.
Recently, Orlik-Rapoport and He
studied this question. Recall that
$x,y\in W$ are {\it $F$-conjugate} if there exists $a \in W$
such that $y=a^{-1} x F(a)$. 
By a case-by-case analysis based on Deligne-Lusztig's criterion,
they obtained the following result (\cite[\S 5]{orlik}, \cite[Theorem 1.3]{he}):

\bigskip

\noindent{\bf Theorem A (Orlik-Rapoport, He).} 
{\it If $w \in W$ is an element of minimal length in its 
$F$-conjugacy class then $\Xb(w)$ is affine.}

\bigskip

When $w$ is a Coxeter element,
the result is due to Lusztig \cite[Corollary 2.8]{Lu2}.
The aim of this note is to generalize, and give a more direct proof
of Theorem A. 
As a consequence (and by applying a combinatorial result on elements 
of minimal length in their $F$-conjugacy class), 
we shall obtain a generalisation of Theorem A. 

\smallskip
Before stating our results, we need some further notation. 
We denote by $B^+$ the {\it braid monoid} associated to $(W,S)$. 
It is the monoid with presentation
$$B^+=\langle (\underline{x})_{x \in W}~|~
\forall~x,x' \in W,~\ell(xx')=\ell(x)+\ell(x') 
\Rightarrow \underline{xx'}=\underline{x}~\underline{x}'\rangle.$$
The automorphism $F$ of $W$ extends to an automorphism of $B^+$ still denoted 
by $F$.

Given $I \subset S$, let $W_I$ denote the subgroup of $W$ generated by $I$ and
let $w_I$ be the longest element of $W_I$ (the element $w_S$ will 
be denoted by $w_0$).
The main result of this note is the following:

\bigskip

\noindent{\bf Theorem B.} 
{\it Let $I$ be an $F$-stable subset of $S$ and 
let $w \in W_I$ be such that there exists a positive integer $d$ and
$a \in B^+$ with
$$\underline{w} F(\underline{w}) \cdots F^{d-1}(\underline{w})=
\underline{w}_I a.$$
Then $\Xb(w)$ is affine.}

\bigskip

The proof of Theorem B is by a general argument while our deduction of
Theorem A relies on combinatorial results on finite Coxeter groups 
(see \cite[Theorem 1.1]{geck michel}, \cite[\S 6]{geck kim} and 
\cite[Theorem 7.5]{he1}) which are proved by a case-by-case analysis. 

\medskip

There is a case where our criterion can be applied easily. 
Indeed, if $d$ is a regular number 
for $(W,F)$ (in the sense of Springer) then by \cite[Proposition 6.5]{BrMi},
there exists a regular element $w \in W$ such that 
$$\underline{w} F(\underline{w}) \cdots F^{d-1}(\underline{w})=\underline{w}_0 
\underline{w}_0.$$
Therefore, by Theorem B, the variety $\Xb(w)$ is affine: this variety is 
of fundamental interest for the geometric version 
of Brou\'e's abelian defect group conjecture for finite reductive 
groups \cite[\S 5.B]{BrMi}.
In particular, if $i \neq j$, this conjecture predicts 
that, as $\qlb\Gb^F$-modules, the cohomology groups 
$H_c^i(\Xb(w),\qlb)$ and 
$H_c^j(\Xb(w),\qlb)$ have no common irreducible 
constituents.

Finally, note that there exists elements satisfying the criterion 
of Theorem B but which do not satisfy Deligne-Lusztig's criterion. 
For instance, if $W$ is of type $B_5$ (and $F$ acts trivially on $W$), 
the element $w=s_1 t s_3 s_2 s_1 t s_1 s_4 s_3 s_2 s_1 t s_1 s_2 s_3$ 
does not satisfy Deligne-Lusztig's criterion (for $q=2$) but satisfies 
$(\underline{w})^5=(\underline{w}_0)^3$ (here, $S=\{t,s_1,s_2,s_3,s_4\}$, 
$ts_1$ has order $4$ and $s_i s_{i+1}$ has order $3$ for $i=1$, $2$, $3$). 
However, this element $w$ is $F$-conjugate by cyclic shift (see Section 2 for 
the definition) to 
$s_4 w s_4=s_1 t s_3 s_2 s_1 t s_1 s_2 s_3 s_4 s_3 s_2 s_1 t s_1$ 
which satisfies Deligne-Lusztig's criterion, 
so the affineness of the variety $\Xb(w)$ can also be obtained 
from Deligne-Lusztig's criterion (see Proposition 2). These computations 
have been checked using {\tt GAP3/CHEVIE} programs written by Jean Michel. 

\medskip
We thank Christian Kaiser and the referee 
for having pointed out a mistake in a previous
version of the paper. We thank Jean Michel for the useful discussions 
we had with him on this subject and for the software programs he 
has written for checking Deligne-Lusztig's criterion.

\bigskip

\section{Preliminaries}

\noindent{\bf Levi subgroup.} 
Let us fix an $F$-stable Borel subgroup $\Bb_0$ of $\Gb$ and
an $F$-stable maximal torus $\Tb$ of $\Bb_0$. Let $\Ub$ be the unipotent
radical of $\Bb_0$.
We identify $N_\Gb(\Tb)/\Tb$ with $W$ by requiring that
$(\Bb_0, w\Bb_0 w^{-1})\in\OCB(w)$. 

Let $I$ be an $F$-stable subset of $S$,
let $\Pb_I=\Bb W_I \Pb$, let $\Vb_I$ denote the unipotent radical of $\Pb_I$ 
and let $\Lb_I$ denote the unique Levi subgroup of $\Pb_I$ containing 
$\Tb$.
Given $w \in W_I$, there
is an isomorphism \cite[Lemma 3]{Lu}
$$\Xb_\Gb(w) \xrightarrow{\sim}
\Gb^F/\Vb_I^F \times_{\Lb_I^F} \Xb_{\Lb_I}(w).$$
In particular,
$$\text{\it $\Xb_\Gb(w)$ is affine if and only if $\Xb_{\Lb_I}(w)$ is affine.}
\leqno{\boldsymbol{(1)}}$$

\bigskip

\noindent{\bf Cyclic shift.} 
If $w$, $w' \in W$, we say that $w$ and $w'$ are 
{\it $F$-conjugate by cyclic shift} 
(and we write $w \cyclic w'$) 
if there exists three sequences $(x_i)_{1 \le i \le n}$, $(y_i)_{1 \le i \le n}$ 
and $(w_i)_{1 \le i \le n+1}$ of elements of $W$ such that
\begin{quotation}
\begin{itemize}
 \item[(1)] $w_1=w$ and $w_{n+1}=w'$;

 \item[(2)] for all $i \in \{1,2,\dots,n\}$, 
$w_i=x_i y_i$, $w_{i+1}=y_i F(x_i)$ and 
$\ell(w_i)=\ell(w_{i+1})=\ell(x_i)+\ell(y_i)$.
\end{itemize}
\end{quotation}
The relation $\cyclic$ is an equivalence relation. 
Two elements which are $F$-conjugate by cyclic shift have the same length. 

\bigskip

\noindent{\bf Proposition 2.} 
{\it If $w \cyclic w'$, then $\Xb(w)$ is affine if and only if 
$\Xb(w')$ is affine.}

\bigskip

\noindent{\sc{Proof}~-}  
By induction, we may assume that there exists $x$ and $y$ in $W$ such that 
$w=xy$, $w'=yF(x)$ and $\ell(w)=\ell(w')=\ell(x)+\ell(y)$. The result follows
from the existence of a purely inseparable morphism $\Xb(w)\to\Xb(w')$
\cite[Page 108]{delu}.\fin

\section{Proof of Theorem B}

Let $I$ be an $F$-stable subset of $S$, let $w \in W_I$ and assume 
that there exists $a \in B^+$ and a positive integer $d$ such that 
$$\underline{w} F(\underline{w}) \cdots F^{d-1}(\underline{w}) = 
\underline{w}_I a.$$
The aim of this section is to prove that $\Xb(w)$ is affine
(Theorem B).
By (1), we may (and we will) assume that 
$I=S$. 

\bigskip

\noindent{\bf Sequences of elements of ${\boldsymbol{W}}$.} 
Given $(x_1,\dots,x_r)$ a sequence of elements of $W$, we set 
$$\OCB(x_1,\dots,x_r) = \OCB(x_1)\times_{\BCB}\cdots\times_{\BCB}
\OCB(x_r).$$

If $(y_1,\dots,y_s)$ is a
sequence of elements of $W$ such that $\underline{x}_1 \cdots \underline{x}_r 
= \underline{y}_1 \cdots \underline{y}_s$ in $B^+$, then 
$\OCB(x_1,\dots,x_r)\simeq \OCB(y_1,\dots,y_s)$ (the varieties are
actually canonically isomorphic \cite[Application 2]{deligne}). For 
a general treatment of these varieties $\OCB(x_1,\dots,x_r)$ 
(and the corresponding Deligne-Lusztig varieties), the 
reader may refer to \cite{DMR}.

\bigskip

\noindent{\bf Proposition 3.} 
{\it If there exists $v \in B^+$ such that 
$\underline{x}_1 \cdots \underline{x}_r = \underline{w}_0 v$, then 
the variety $\OCB(x_1,\dots,x_r)$ is affine.}

\bigskip

\noindent{\sc{Proof}~-} 
Let $v_1$,\ldots, $v_n\in W$ be such that 
$\underline{v}_1 \cdots  \underline{v}_n = v$. We have 
$\OCB(x_1,\dots,x_r) \simeq \OCB(w_0,v_1,\dots,v_n)$,
so it remains to prove that $\OCB(w_0,v_1,\dots,v_n)$ is 
affine. 

For each $x \in W$, we fix a representative $\xdo$ of $x$ in $N_\Gb(\Tb)$. 
We set 
$$\tilde{\OCB}(x_1,\dots,x_r) = 
\{(g_0\Ub,g_1\Ub,\dots,g_r\Ub) \in (\Gb/\Ub)^{r+1}~|~
\forall~1 \le i \le r,~g_{i-1}^{-1} g_i \in \Ub \xdo_i \Ub\}.$$

The group $\Tb$ 
acts on the right on $\tilde{\OCB}(x_1,\dots,x_r)$ as follows: 
$$(g_0\Ub,g_1\Ub,\dots,g_r\Ub) * t = 
(g_0t\Ub,g_1 \lexp{x_1^{-1}}{t}\Ub,\dots,
g_r\lexp{x_r^{-1}\cdots x_1^{-1}}{t} \Ub).$$

The canonical map 
$$\fonctio{\tilde{\OCB}(x_1,\dots,x_r)}{\OCB(x_1,\dots,x_r)}{
(g_0\Ub,g_1\Ub,\dots,g_r\Ub)}{(g_0\Bb_0g_0^{-1},g_1\Bb_0g_1^{-1},
\dots,g_r\Bb_0g_r^{-1})}$$
identifies $\OCB(x_1,\dots,x_r)$ with the quotient of 
$\tilde{\OCB}(x_1,\dots,x_r)$ 
by $\Tb$: indeed, both varieties are smooth (hence normal), 
the above map is smooth (hence separable) and it is easily checked 
that its fibers are precisely the $\Tb$-orbits. 
Since $\Tb$ acts freely on $\tilde{\OCB}(x_1,\dots,x_r)$, and since 
the quotient of an affine variety by a free action of a torus is affine,
\cite[Corollary 8.21]{borel}, 
the result will follow if we are able to prove that
$\tilde{\OCB}(w_0,v_1,\dots,v_n)$ is affine. Therefore, it is sufficient 
to show that the map 
$$\begin{array}{ccl}
\varphi:\Gb \times \displaystyle{\prod_{i=1}^n} 
(\Ub \vdo_i \cap \vdo_i \Ub^-)\!\!\! & \longto 
& \!\!\!\tilde{\OCB}(w_0,v_1,\dots,v_n)\\
(g;h_1,\dots,h_n) & \longmapsto & 
\!\!\!(g\Ub,g\wdo_0\Ub,g\wdo_0 h_1\Ub,g\wdo_0 h_1 h_2\Ub,\dots,
g\wdo_0 h_1\cdots h_n\Ub)
  \end{array}$$
is an isomorphism of varieties (here, $\Ub^-=\lexp{w_0}{\Ub}$). 

In order to prove that $\varphi$ is an isomorphism, 
we shall construct its inverse. For this, we shall need some notation. 
First, the map $\Ub \times \Ub \to \Ub \wdo_0 \Ub$, 
$(x,y) \mapsto x \wdo_0 y$ is an isomorphism of varieties: 
we shall denote by $\Ub \wdo_0 \Ub \to \Ub \times \Ub$, $g \mapsto 
(\eta(g),\eta'(g))$ its inverse. Also, the map 
$\Ub \vdo_i \cap \vdo_i \Ub^- \times \Ub \to \Ub \vdo_i \Ub$, $(x,y) 
\mapsto xy$ is an isomorphism of varieties ($i=1$, $2$, \dots, $n$), 
and we shall denote by 
$\eta_i : \Ub \vdo_i \Ub \to \Ub \vdo_i \cap \vdo_i \Ub^-$ 
the composition ot its inverse with the first projection. 
Note that, if $g \in \Ub \wdo_0 \Ub$, $h \in \Ub \vdo_i \Ub$ and 
$u$, $v \in \Ub$, then
$$\eta(ugv)=u\eta(g),\quad\eta(g)\wdo_0\Ub=g\Ub,\quad\eta_i(hv)=\eta_i(h)
\quad\text{and}\quad\eta_i(h)\Ub=h\Ub.\leqno{(*)}$$
Now, if $x=(g\Ub,g_0\Ub,g_1\Ub,\dots,g_n\Ub) \in
\tilde{\OCB}(w_0,v_1,\dots,v_n)$, we set 
\begin{eqnarray*}
\psi(x)&=& g \eta(g^{-1}g_0),\\
\psi_0(x)&=& \psi(x) \wdo_0, \\
\psi_i(x)&=& \eta_i\bigl((\psi_0(x)\psi_1(x)\cdots\psi_{i-1}(x))^{-1} g_i\bigr),
\end{eqnarray*}
for all $i \in \{1,2,\dots,n\}$. By $(*)$, 
the maps $\psi$ and $\psi_j$ are well-defined morphisms 
of varieties and it is easily checked that the morphism of varieties 
$$\begin{array}{ccl}
\tilde{\OCB}(w_0,v_1,\dots,v_n) & \longto 
& \Gb \times \displaystyle{\prod_{i=1}^n} (\Ub \vdo_i \cap \vdo_i \Ub^-)\\
x & \longmapsto & (\psi(x);\psi_1(x),\dots,\psi_n(x))
\end{array}$$
is well-defined and is an inverse of $\varphi$.\fin

\bigskip

\noindent{\bf Introducing Frobenius.} 
The morphism
$$\Xb(w)\to \BCB^d,\ \Bb\mapsto
(\Bb,F(\Bb),\ldots,F^{d-1}(\Bb))$$
indentifies $\Xb(w)$ with the closed subvariety 
$\Delta_d\cap \OCB(w,F(w),\dots,F^{d-1}(w))$,
where
$\Delta_d=\{(\Bb,F(\Bb),\dots,F^{d-1}(\Bb))~|~\Bb\in\BCB\}$ is
a closed subvariety of $\BCB^d$.
By Proposition 3,
the variety $\OCB(w,F(w),\dots,F^{d-1}(w))$ is affine, hence
$\Xb(w)$ is affine as well. 
The proof of Theorem B is complete. 

\section{Proof of Theorem A}

Let $C$ be an $F$-conjugacy class in $W$ and $C_{\min}$
its subset of elements of minimal length.
Let $d$ be the smallest positive integer $k$ such that
$w F(w) \dots F^{k-1}(w)=1$ and $F^k$ acts as the identity 
on $W$ for $w\in C_{\min}$.
Following \cite[Theorem 1.1]{geck michel} (in the split case) 
and \cite[Definition 5.3]{geck kim} (in the general case), we say that
an element $w \in C_{\min}$ 
is {\it good} if there exists a sequence 
$I_1 \supseteq I_2 \supseteq \cdots \supseteq I_r$ of subsets 
of $S$ such that 
$$\underline{w} F(\underline{w}) \dots F^{d-1}(\underline{w}) = 
\underline{w}_{I_1}^2 \underline{w}_{I_2}^2 
\cdots \underline{w}_{I_r}^2\leqno{(*)}$$
in $B^+$. 

\bigskip

\noindent{\bf Proposition 4.} 
{\it If $w$ is a good element of $C_{\min}$, then $\Xb(w)$ is affine.} 

\bigskip

\noindent{\sc{Proof}~-}  By Theorem B, it remains to show that the subset $I_1$ 
of the identity $(*)$ is $F$-stable. Let $I$ be the set of simple 
reflections occuring in a reduced expression of $w$ (note that 
$I$ does not depend on the choice of the reduced expression 
\cite[Corollary 1.2.3]{geck pfeiffer}). Then the set of $s \in S$ such that 
$\underline{s}$ occurs in a reduced expression of 
$\underline{w} F(\underline{w}) \dots F^{d-1}(\underline{w})$ 
is equal to $I \cup F(I) \cup \cdots \cup F^{d-1}(I)$ 
(by looking at the left-hand side of $(*)$) and is also equal 
to $I_1$ (by looking at the right-hand side). Since $F^d$ acts 
as the identity on $W$, we get that $I_1$ is $F$-stable.\fin 

\bigskip

Let us now come back to the proof of Theorem A. Let $w \in C_{\min}$. 
Let $I$ be the minimal $F$-stable subset of $S$ such that $w \in W_I$. 
By (1), we may assume that $I=S$. Now, if $w' \in C_{\min}$, 
then $w \cyclic w'$ (see \cite[Theorem 3.2.7]{geck pfeiffer} for the 
split case, \cite[\S 6]{geck kim} for twisted exceptional groups and 
\cite[Theorem 7.5]{he1} for twisted classical groups), so 
$\Xb(w)$ is affine if and only if $\Xb(w')$ is affine 
by Proposition 2. 
Therefore, the result follows from Proposition 4 and the next Theorem:

\bigskip

\noindent{\bf Theorem 6 (Geck-Michel, Geck-Kim-Pfeiffer, He).} 
{\it There exists a good element in $C_{\min}$.}

\bigskip

\noindent{\sc Proof - } 
By standard arguments (see \cite[\S 5.5]{geck kim}), we 
may assume that $W$ is irreducible. If $F$ acts trivially 
on $W$, the Theorem is \cite[Theorem 1.1]{geck michel}. 
If $F$ does not act trivially and $W$ is not of type $A$, this is
\cite[\S 5.5]{geck kim}. When $W$ is of type $A$ and 
$F$ acts non trivally on $W$, this follows from
\cite[Corollary 7.25]{he1}.\fin

\end{document}